\documentclass[conference]{IEEEtran}
\IEEEoverridecommandlockouts
\usepackage{cite}
\usepackage{amsmath,amssymb,amsfonts}
\usepackage{algorithmic}
\usepackage{graphicx}
\usepackage{textcomp}
\usepackage[english]{babel}
\usepackage[utf8]{inputenc}
\usepackage{multirow}
\usepackage{tabularx}
\usepackage{nomencl}
\makenomenclature
\def\BibTeX{{\rm B\kern-.05em{\sc i\kern-.025em b}\kern-.08em
    T\kern-.1667em\lower.7ex\hbox{E}\kern-.125emX}}
\usepackage{etoolbox}
\renewcommand\nomgroup[1]{%
  \item[\bfseries
  \ifstrequal{#1}{V}{Variables}{%
  \ifstrequal{#1}{P}{Parameters}{%
  \ifstrequal{#1}{S}{Sets}{}}}%
]}

\title{Modeling A Micro-Nexus of Water and Energy\\for Smart Villages/Cities/Buildings}

\begin{document}
\author{

\IEEEauthorblockN{Qifeng Li$^\ast$,~\IEEEmembership{Member,~IEEE,} Suhyoun Yu$^\ast$, Ameena Al-Sumaiti$^\dagger$, and Konstantin Turitsyn$^\ast$,~\IEEEmembership{Member,~IEEE,}}
\IEEEauthorblockA{$^\ast$Department of Mechanical Engineering, Massachusetts Institute of Technology, Cambridge, MA USA \\
(e-mail: \{qifengli, syu2, turitsyn\}@mit.edu)\\
$^\dagger$Dept. of Electrical \& Computer Engr., Masdar Inst., Khalifa Univ. Sci. \& Tech., Abu Dhabi, UAE \\
(e-mail: aalsumaiti@masdar.ac.ae)}

\thanks{This work is supported by the MI-MIT Cooperative Program under grant MM2017-000002.}
}

\maketitle

\begin{abstract}
This paper introduces a micro-nexus of water and energy which can be considered as one of the physical infrastructures of the future building/city/village systems. For the electricity side, an alternating current (AC) power flow model integrated with battery energy storage and renewable generation is adopted. The nonlinear hydraulic characteristics in pipe networks is also considered in the proposed micro water-energy nexus (WEN) model. Integer variables are involved to represent the on/off state of pumps. Base on the proposed nexus model, a co-optimization framework of water and energy networks is developed. The overall co-optimization model is a mixed-integer nonlinear programming problem which is tested on a water-energy nexus which consists of the IEEE 13-bus distribution system and a 8-node water distribution network. The simulation results demonstrate that the cost-efficiency of the co-optimization framework is higher than optimizing two systems independently.
\end{abstract}


\section{Introduction}

The smart building/city/village/ are development visions of buildings/cities/villages of which an essential function is to improve the efficiency and security of various services by integrating some emerging technologies (such as Internet of Things, optimization, information, communication, and control) \cite{Albino}- \cite{Ameena}. In a smart building, city or village, many physical systems are interconnected, which enables co-operation or co-optimization of these systems. The water and energy networks are two of the most important physical systems in almost all types of communities, including cities, remote villages, as well as buildings, since efficient and secure water and energy services are considered as life necessities. Therefore, it is reasonable to consider an intelligent infrastructure of water and energy as one of the foundations of the smart buildings/cities/villages. In fact, refs. \cite{Morvaj} and \cite{Heap} have demonstrated the importance of a smart energy grid for the smart building/city/village. 

It is no coincidence that water and energy systems are tightly intertwined. On one hand, most of the services provided by the water system consume energy. Water supply, seawater desalination, groundwater pumping, and wastewater treatment even account for a considerable amount of total energy consumption in some area like Middle-East and North Africa (MENA). Many countries in MENA consume 5-12$\%$ or more of total electricity consumption for desalination \cite{Siddiqi}. On the other hand, the intermittent renewable energy resources and electro-mobility introduce critical uncertainty in achieving the balance of supply and demand in the power system \cite{Palensky}. The services of water systems can be considered as a source of grid flexibility to harness such uncertainty \cite{Menke}. The energy-water nexus has attracted substantial research efforts recently as more and more researchers realized the mutual influence of water and energy systems \cite{Menke} - \cite{Zhang}.

Demand side management (DSM) offers promising measures to mitigate the imbalance between electric supply and demand by utilizing the flexibility of loads. To explore the potential of water distribution systems (WDS) in providing demand response capability, we introduce a micro-nexus of water and energy on the demand side. A new mathematical model needs to be developed in order to better capture the interactions between water and electric distribution systems. The optimal pump scheduling problem (OPS) has been studied for decades in the water sector \cite{Jowitt}. In this problem, researchers and engineers pursue an operation schedule of pumps that can satisfy the requirement of water supply with the minimum energy consumption. No electricity networks are considered in this problem. Recently, the interactions between the water and energy systems start receiving researchers' attentions. Ref. \cite{Santhosh} defines a supply side water-energy nexus (WEN) assuming that water transmission systems couple with power transmission grids only at the co-generation plants. The DC power flow model and linear pipe network model considered in \cite{Santhosh} are not valid for the distribution-level nexus. A recent paper \cite{Menke} explores the opportunities for WDS to provide frequency response services through pump scheduling. In essence, the approach proposed in \cite{Menke} is still an optimal pump scheduling problem without considering the electricity network. 

The micro water-energy nexus introduced in this paper comprises microgrids and micro-WDSs which are directly connected with customers. Nonlinear network models with high-fidelity are considered for both water and power systems. Such a nexus model lies a solid foundation for many future research interests such as optimal water-DSM schemes, co-optimization of water and energy systems, and co-security/reliability of micro-WEN. It is also beneficial for developing a coordinated islanding scheme of electric and water networks under emergencies.

\section{Design of the Micro Water-Energy Nexus}

\subsection{Background of Water Distribution System}

Water distribution systems deliver water to consumers with  acceptable quality, appropriate pressure, and the required quantity, which are considered as the backbone of the society since water is one of the prime needs of humans. Within the next decade, approximately 1.8 billion people worldwide will be living in areas of absolute water scarcity. Such a situation requires that the function of a future WDS is not just the distribution of water through pipe networks. It involves water treatment, water recycling, water purification, cooling, metering/monitoring, funds/revenue generated, etc \cite{Kashid}. That means the water and energy systems will be more strongly intertwined in the future. The integration of intelligent metering systems generates the ideal of "intelligent" or "smart" WDS \cite{Mizuki}. It enables the strategic co-operation of all facilities in a WDS. SWAN\footnote{SWAN: The Smart Water Networks Forum (www.swan-forum.com).} defines a smart water network as a collection of data-driven components helping to operate the data-less physical layer of pipes, pumps, reservoirs and valves. Water utilities are gradually deploying more data-enabled components. It is up to us to make the most out of them, by turning the discrete elements into a cohesive "overlay network".

\subsection{Micro Water-Energy Nexus}

\begin{figure}[tb]
\centering
\includegraphics[width=0.4\textwidth]{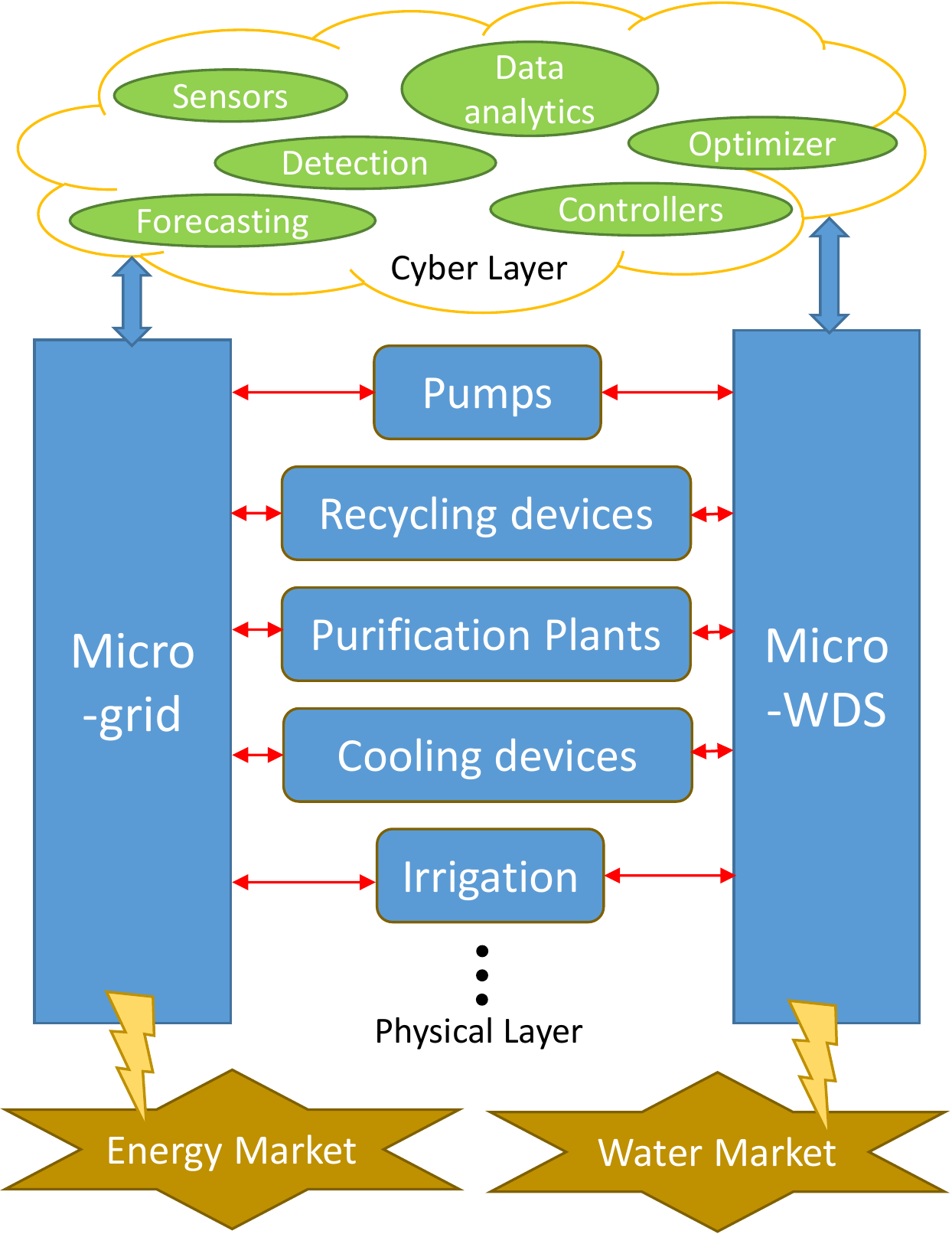}
  \caption{Schematic of a smart micro-WEN.}
  \label{fig:cyberphysical}
  
\end{figure}
\begin{figure*}[tb]
\centering
\includegraphics[width=0.8\textwidth]{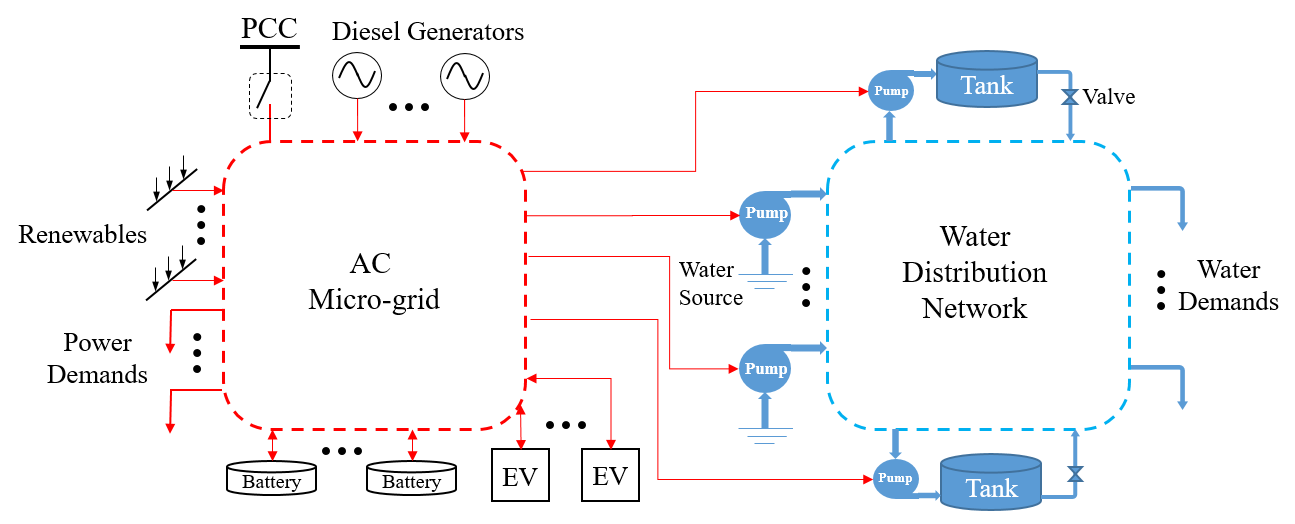}
  \caption{A physical structure of the micro water-energy nexus.}
  \label{fig:WEN}
\end{figure*}

Similar to the microgrid, a micro-WDS is also a cyber-physical system. Therefore, it is reasonable to design a micro-nexus of water and energy in which the physical layers of the two systems are interconnected, so that they could share the main parts of the cyber layers as shown in Figure \ref{fig:cyberphysical}. Microgrid technologies, like battery energy storage systems (BESS), high-inertia synchronous generators, various power electronic devices, and electric vehicles (EV), are also considered in the micro-WEN which is applicable to many types of communities.

\textbf{Remote Villages}. The micro-WEN of a smart village is generally off-grid and its energy consumption mainly comes from renewable generation. High-inertia synchronous generators like diesel generators serve to maintain the stability of the system frequency rather than provide electric energy. The water source may be wells, or sea water desalination plants if the remote village is located on the coast.

\textbf{Blocks in cities}. In the development vision of cities, renewable energy accounts for a considerable portion of energy consumed. The micro-WEN of a block in a city is connected to the power grid via the point of common coupling (PCC) under normal condition and operated off-grid when the service of main power system is not available. The water source is primarily the higher level water supply system. The micro-WEN needs to have sufficient capability to support EV chargers no matter in the grid-connected mode or the islanding mode.

\textbf{Buildings}. The case in buildings is similar to that of blocks in cities. The Physical structure of the water network in a building is different from those in a village or a city block. However, they share the similar mathematical model.

\section{Modeling of the Micro Water-Energy Nexus}

\nomenclature[V]{$P_{ij,t}$}{Active branch power in line $ij$}
\nomenclature[V]{$Q_{ij,t}$}{Reactive branch power in line $ij$}
\nomenclature[V]{$e_{i,t}$}{Real part of voltage at bus $i$}
\nomenclature[V]{$f_{i,t}$}{Imaginary part of voltage at bus $i$}
\nomenclature[V]{$P_{i,t}^G $}{Active power of diesel generator at bus $i$}
\nomenclature[V]{$Q_{i,t}^G $}{Reactive power of diesel generator at bus $i$}
\nomenclature[V]{$P_{i,t}^{ES} $}{Active power of the BESS at bus $i$}
\nomenclature[V]{$ Q_{i,t}^{ES}$}{Reactive power of the BESS at bus $i$}
\nomenclature[V]{$P_{i,t}^{Pump} $}{ Active consumptions of the pump at bus $i$}
\nomenclature[V]{$ Q_{i,t}^{Pump}$}{Reactive consumptions of the pump at bus $i$}
\nomenclature[V]{$L_{i,t}^{ES} $}{ Active power loss of the BESS at bus $i$}
\nomenclature[V]{$x_t $}{ Vector of water flow in pipes}
\nomenclature[V]{$y_t $}{ Vector of water head at nodes}
\nomenclature[V]{$w_{i,t}^S $}{ Water flow to the water tank at node $i$}
\nomenclature[V]{$w_{i,t}^G $}{ Water flow injected from the water source in pipe $i$}
\nomenclature[V]{$y_{ij,t}^G $}{ Head gain imposed by the pump in pipe $ij$}
\nomenclature[V]{$\alpha_{i,t}$}{Binary variable denoting the on/off state of the pump in pipe $i$}
\nomenclature[P]{$A $}{Incidence matrix}
\nomenclature[P]{$a_{ij}$}{Coefficient of pump characteristics in pipe $ij$}
\nomenclature[P]{$b_{ij}$}{Coefficient of pump characteristics in pipe $ij$}
\nomenclature[P]{$P_{i,t}^L $}{Active electric load at bus $i$}
\nomenclature[P]{$Q_{i,t}^L$}{Reactive electric load at bus $i$}
\nomenclature[P]{$P_{i,t}^{RE} $}{Active power output of renewable at bus $i$}
\nomenclature[P]{$Q_{i,t}^{RE}$}{Reactive power output of renewable at bus $i$}
\nomenclature[P]{$G_{ij}$}{ Conductance of line $ij$}
\nomenclature[P]{$B_{ij}$}{ Susceptance of line $ij$}
\nomenclature[P]{$r_i^{Batt}$}{ Loss coefficient related to the battery of BESS}
\nomenclature[P]{$r_i^{Cvt}$}{ Loss coefficient related to the converter of BESS}
\nomenclature[P]{$E_{i,0}^{ES}$}{ Initial state of charging of BESS at bus $i$}
\nomenclature[P]{$h_i $}{Elevation at node $i$ of the water network}
\nomenclature[P]{$R_{ij}^P$}{ Head loss coefficient of pipe $ij$}
\nomenclature[P]{$c_t$}{Locational marginal price at PCC}
\nomenclature[P]{$c_{1,i}$}{Coefficient of the first-order term in the cost function of diesel generator at bus $i$}
\nomenclature[P]{$c_{2,i}$}{Coefficient of the second-order term in the cost function of diesel generator at bus $i$}
\nomenclature[S]{$\mathcal{E}_E$}{Edge set of the electricity network}
\nomenclature[S]{$ \mathcal{E}_W$}{Edge set of the water network}
\nomenclature[S]{$\mathcal{N}_E$}{Bus set of the electricity network}
\nomenclature[S]{$ \mathcal{N}_W$}{Node set of the water network}
\nomenclature[S]{$\mathcal{E}_W^P$}{Set of pipes with a pump installed}
\nomenclature[S]{$\mathcal{N}_E^S$}{Set of bus with a BESS connected}
\nomenclature[S]{$ \mathcal{N}_W^S$}{Set of nodes connected to a tank}
\nomenclature[S]{$\mathcal{N}_E^G$}{Set of bus with controllable generations}
\nomenclature[S]{$\mathcal{N}_E^P$}{Set of bus with a pump connected}

\printnomenclature

This section develops a mathematical model of the introduced physical layer of micro-WEN as shown in Figure \ref{fig:WEN}. The problem considered in this paper is of multi-period. Unless otherwise stated, the subscript $t$ denotes the time period. For the AC-microgrid, we use the quadratic AC power flow model \cite{LiThe} with renewable generation and batteries integrated to describe the electricity network. The detailed model is given as
\begin{subequations}
\begin{gather} 
P_{i,t}^G+P_{i,t}^{RE}-P_{i,t}^{Pump}-P_{i,t}^L-P_{i,t}^{ES} =\sum_jP_{ij,t} \\
Q_{i,t}^G+Q_{i,t}^{RE}-Q_{i,t}^{Pump}-Q_{i,t}^L-Q_{i,t}^{ES}= \sum_jQ_{ij,t} \\
P_{ij,t}=G_{ij}(e_{i,t}e_{j,t}+f_{i,t}f_{j,t})-B_{ij}(e_{i,t}f_{j,t}+f_{i,t}e_{j,t})\\
Q_{ij,t}=G_{ij}(f_{i,t}e_{j,t}-e_{i,t}f_{j,t})-B_{ij}(e_{i,t}e_{j,t}+f_{i,t}f_{j,t}) \\
P_{ij,t}^2+Q_{ij,t}^2 \le \overline{S}_{ij}^2 \\
\underline{V}^2_i \le e_{i,t}^2+f_{i,t}^2 \le \overline{V}_i^2 \\
\underline{P}_k^G \le P_{k,t}^G \le \overline{P}_k^G \\
\underline{Q}_k^G \le Q_{k,t}^G \le \overline{Q}_k^G,
\end{gather}
\end{subequations}
where $i,j \in \mathcal{N}_E$, $ij \in \mathcal{E}_E$, and $k \in \mathcal{N}_E^G$. Constraints (1a) and (1b) describe the nodal balance of active and reactive power respective; (1c) and (1d) stem from the Ohm's law; (1e)-(1h) are system constraints.

The following nonlinear model of battery energy storage unit is adopted for the mathematical model of micro-WEN. Please refer to \cite{LiConvex} for more details about this BESS model. For $\forall i \in \mathcal{N}_E^S$, we have
\begin{subequations}
\begin{gather} 
(r_i^{Batt} + r_i^{Cvt})(P_{i,t}^{ES})^2 + r_i^{Cvt}(Q_{i,t}^{ES})^2 = L_{i,t}^{ES}(e_{i,t}^2+f_{i,t}^2) \\
(P_{i,t}^{ES})^2 + (Q_{i,t}^{ES})^2 \leqslant (\overline{S}_{i}^{ES})^2 \\
\underline{E}_{i}^{ES} \leqslant E_{i,0}^{ES} - \int^t_0 (P_{i,t}^{ES}+L_{i,t}^{ES})\,d\tau \leqslant \overline{E}_{i}^{ES}.
\end{gather}
\end{subequations}

The micro-WDS in the designed micro-WEN comprises a pipe network, pumps, water storage (i.e. tanks), and loads. For the micro-WDS, we have the following assumptions:\\
1) \textit{The pipe network is a directed graph $\mathcal{G}_W=(\mathcal{N}_W,\mathcal{E}_W)$ and, hence, $A_{ij} \in \{-1,\ 0,\ 1\}$ for all pipes};\\
2) \textit{A pump is located in one of the pipes and it imposes a head gain to the pipe when it is on or, otherwise, the pipe is considered to be closed};\\
3) \textit{The pumps are constant-speed motors and the water head gain introduced by the pump is a linear function of the water flow goes through the pump which is given by}\footnote{The hydraulic characteristics of a constant-speed pump is generally approximated by a quadratic function of the water flow across the pump, i.e. $y=cx^2+ax+b$. Coefficient $c$ is usually very small compared with $a$ and $b$ \cite{UlanickiModeling}. Hence, we simply make $c=0$ in this paper.}
\begin{equation}
y_{ij,t}^G = a_{ij}x_{ij,t}+b_{ij} \quad ((i,j) \in \mathcal{E}_W^P);
\end{equation}
4) \textit{The pump converts electric power into mechanical power at a constant efficiency of} $\eta$;\\
5) \textit{The power factors of pumps are fixed, namely} $P_{ij,t}^{Pump}/Q_{ij,t}^{Pump}$=\textit{constant}.

Under the above assumptions, the mathematical model of the micro-WDS can be expressed as
\begin{subequations}
\begin{gather} 
\sum_{(i,j) \in \mathcal{E}_W}A_{k,ij}x_{ij,t} = w_{k,t}^G - w_{k,t}^S - w_{k,t}^D,(k \in \mathcal{N}_W) \\
y_{i,t}-y_{j,t}+h_i-h_j= R_{ij}^Psgn(x_{ij,t})x_{ij,t}^2,((i,j) \in \mathcal{E}_W/\mathcal{E}_W^P)\\
\begin{cases}
\begin{array}{ll}
 \begin{split}
 &y_{i,t}-y_{j,t}+h_i-h_j\\&+y_{ij,t}^G= R_{ij}^Px_{ij,t}^2
 \end{split} & \text{if} \ \alpha_{ij,t}=1 \\
 x_{ij,t} = 0, & \text{if} \ \alpha_{ij,t}=0
\end{array}, ((i,j) \in \mathcal{E}_W^P) 
\end{cases} \\
 \underline{S}_{i}^{w} \leqslant S_{i,0}^w + \int^t_0 w_{i,\tau}^S \,d\tau \leqslant \bar{S}_{i}^{w},\quad (i \in \mathcal{N}_W^S)\\
 \underline{x} \le x_t \le \overline{x}, \\
 \underline{y} \le y_t \le \overline{y}, \\
 \underline{w}_i^G \le w_{i,t}^G \le \overline{w}_i^G, \  (i \in \mathcal{N}_W^G) \\
 \underline{w}_i^S \le w_{i,t}^S \le \overline{w}_i^S, \  (i \in \mathcal{N}_W^S)
\end{gather}
\end{subequations}
where $i,j \in \mathcal{N}_W$; equation (4a) represents the mass balance of the water network; constraints (4b) and (4c) formulate the hydraulic characteristics of a normal pipe and the pipe with a pump installed respectively; constraint (4d) denotes the state of charing of the water tanks; (4e)-(4h) are system constraints; $syn(x)=-1$ if $x \le 0$ or, otherwise $syn(x)=1$. 

The following constraints act as the mathematical link between the microgrid ((1)-(2)) and the micro-WDS ((3)-(4))
\begin{equation}
\eta P_{k,t}^{Pump}=x_{ij,t}y_{ij,t}^G=a_{ij}x_{ij,t}^2+b_{ij}x_{ij,t},
\end{equation}
where $k \in \mathcal{N}_E^P$ and $(i,j) \in \mathcal{E}_W^P$.

\textbf{Remark:} (1)-(5) together form the mathematical model of the designed micro-WEN. The model provides a framework in which the stochasticity of renewable generation, water and electricity demands can be easily incorporated whenever needed. It will not be difficult to incorporate more electric and water appliances, such as electric vehicles and water purification, irrigation, desalination, recycling devices, into the developed mathematical model.

\section{Co-optimization of Water and Energy Networks}
Based on the mathematical model of micro-WEN introduced in Section III, this section develops a co-optimization framework for water and energy networks. The objective of this co-optimization problem is to minimize the total energy cost for meeting the demands of both electricity and water. We formulate the energy cost as

\begin{equation}
f(P_{i,t}^G)=\sum_t(c_tP_{1,t}^G + \sum_{i \in \mathcal{N}_E^G/\text{PCC}}(c_{1,i}P_{i,t}^G+c_{2,i}(P_{i,t}^G)^2)),
\end{equation}
where $P_{1,t}^G$ denotes the power from the grid via PCC (i.e. the serial number of PCC is 1); $c_t$ can be consider as the nodal prices at PCC which are calculated as the result of security constrained economic dispatch (SCED) by ISOs/RTOs. As a result, the co-optimization model is
\begin{align}
\begin{split}
\min_{P_{i,t}^G} \quad &\text{(6)} \\
\text{s.t.} \quad &\text{(1)-(5)}
\end{split}. \tag{CO-OPT}
\end{align}

\section{Case Study}

To evaluate the efficiency of the proposed co-optimization framework (CO-OPT), we compare it with an existing daily operation scheme. We assume that this scheme comprises two main stages. In the first stage, the water utility tries to minimize the energy consumption by doing a day-ahead OPS based on the day-ahead water demand forecast, which is given by

\begin{align}
\begin{split}
\min_{\alpha_{i,t}} \quad &\sum_t(\sum_{(i,j) \in \mathcal{E}_W^P}P_{ij,t}^{Pump}) \\
\text{s.t.} \quad &\text{(3) and (4)}
\end{split}. \tag{OPS}
\end{align}

The second stage is an optimization problem which is analogous to the unit commitment (UC) in power transmission systems. The UC-like problem is performed day-ahead by the power utility to minimize the total energy cost for meeting all electricity demands based on the day-ahead forecast of electricity demands, where the diesel generators and BESS units are control devices. With the output of (OPS), i.e. $P_{ij,t}^{Pump}\ ((i,j) \in \mathcal{E}_W^P)$, as a set of time-varying parameters $P_{k,t}^{Pump}\ (k \in \mathcal{N}_E^P)$, the second stage problem is
\begin{align}
\begin{split}
\min_{P_{i,t}^G} \quad &\text{(6)} \\
\text{s.t.} \quad &\text{(1) and (2)}
\end{split}. \tag{UC}
\end{align}

\begin{figure}[t]
\centering
\includegraphics[width=0.5\textwidth]{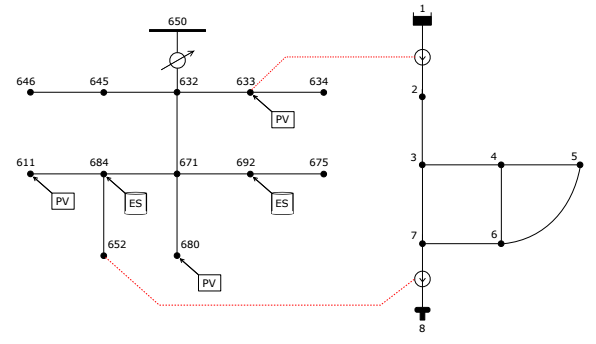}
  \caption{Topology of the test system.}
  \label{fig:Test}
\end{figure}

The micro-WEN for the case study is composed of the IEEE 13-bus system and an 8-node WDS from the EPANET manual \cite{EPANET}. The topology of the test micro-WEN is given in Figure \ref{fig:Test}. We assume that the 13-bus microgrid is integrated with high penetration of PV resources. The pumps deliver 45.72 meter of head at a flow of 0.038 m$^3$/s. The tank is 18.3 meters in diameter and 5.1 meters in depth. The 24-hour load profile of the 13-bus system is generated by applying a typical 24-hour summer loadshape (please refer to Figure 41 in \cite{NarangHigh}). For the 24-hour demand profiles of the water system, please refer to the EPANET manual \cite{EPANET}. Further details about the PV systems and BESSs are given in Table \ref{Table2}.

\begin{table}[h]
\centering
\caption{PV system and BESS Location and Capacity}
\label{Table2}
\begin{tabular}{cc}
\hline\hline
\textbf{PV location (bus \#) and capacity} & \textbf{Penetration} \\
633 (0.5 MW), 680 (0.2 MW), 684 (0.5 MW)   & 36.7\%               \\
\hline
\multicolumn{2}{c}{\textbf{DES location (bus \#) and capacity}}   \\
\hline
\multicolumn{2}{c}{684 (0.6 MW, 2.4 MWh), 692 (0.8 MW, 3.2 MWh)} \\
\hline\hline
\end{tabular}
\end{table}

This section compares the optimal solutions of two optimization problems. Problem (i) is a two stage problem which consists of (OPS) as the first stage problem and (UC) as the one in second stage. Problem (ii) is the proposed framework (CO-OPT). The problems are solved by the corresponding solvers, as shown in Table \ref{table1}, through the optimization package JuMP in Julia (version 0.5.2). It can be observed from the results tabulated in Table \ref{table1} that the co-optimization framework can produce higher cost-efficiency.
\begin{table}[h]
\centering
\caption{Problem Classification, Solvers, and Results}
\label{table1}
\begin{tabular}{ccccc}
\hline\hline
\multicolumn{2}{c}{\textbf{Problem}} & \textbf{Classification} & \textbf{Solver} & \multicolumn{1}{l}{\textbf{\begin{tabular}[c]{@{}l@{}}Optimal\\ Solution\end{tabular}}} \\
\hline
\multirow{2}{*}{(i)}    & (OPS)      & MINLP\textsuperscript{1}   & BONMIN\cite{Grossmann}          & \multirow{2}{*}{1445}                                                                    \\
                        & (UC)       & NLP\textsuperscript{2}     & IPOPT\cite{IPOPT}           &                                                                                         \\
(ii)                    & (CO-OPT)   & MINLP                   & BONMIN          & 1339 \\
\hline\hline
\multicolumn{5}{l}{\textsuperscript{1}\footnotesize{NINLP: Mixed-integer Nonlinear Programming}} \\
\multicolumn{5}{l}{\textsuperscript{2}\footnotesize{NLP: Nonlinear Programming}}
\end{tabular}
\end{table}

\section{Conclusion and Future Work}

The coupling in the demand side of water and energy systems is becoming stronger and stronger as the number of electrified water devices (such as water recycling systems, desalination facilities, and water purification devices) is increasing. This paper designs a micro water-energy nexus which can be considered as the physical foundation of smart buildings/cities/villages. The introduced micro-WEN is composed of a microgrid and a micro water distribution system. A mixed-integer nonlinear mathematical model is developed to formulate the behaviors of the micro-WEN. Currently, the power and water systems are operated independently by different utilities. A numerical study in this paper shows that a higher cost-efficiency can be achieved by operating the water and energy networks as a whole system. The introduced physical and mathematical models of micro-WEN will facilitate many future research which is related to smart cities/villages/buildings. Following this fundamental development, we will continue to explore effective demand response schemes of WDSs by considering more water system services like irrigation. Means of improving the security and reliability for micro-WENs will be investigated by fully considering the interactions between water and electricity systems.

\end{document}